\input amstex
\loadbold
\documentstyle{amsppt}

\topmatter
\title{Universal Schubert Polynomials}
\endtitle
\author{William Fulton}\endauthor
\address{University of Chicago, Chicago IL 60637} \endaddress
\email fulton\@math.uchicago.edu \endemail
\thanks The author was supported by a Tage Erlander Guest 
Professorship at Institut Mittag-Leffler and the National Science 
Foundation 
\abstract
We introduce some polynomials that represent general degeneracy loci 
for maps of vector bundles, and which specialize to the known 
classical and quantum forms of single and double Schubert polynomials.
\endabstract
\endthanks
\date {February 15, 1998} \enddate
\endtopmatter

\document

\head 1. Introduction \endhead

The aim of this paper is to introduce some polynomials that specialize 
to all previously known Schubert polynomials: the classical Schubert 
polynomials of Lascoux and Sch\" utzenberger \cite{L-S1}, \cite{M}, the 
quantum Schubert polynomials of Fomin, Gelfand, and Postnikov 
\cite{F-G-P}, and quantum Schubert polynomials for partial flag 
varieties of Ciocan-Fontanine \cite{CF2}.  There are also double 
versions of these universal Schubert polynomials that generalize the 
previously known double Schubert polynomials \cite{L}, \cite{M}, 
\cite{K-M}, \cite{CF-F}.  They describe degeneracy loci of maps of 
vector bundles, but in a more general setting than the previously known 
setting of \cite{F2}.  

These universal Schubert polynomials possess many but not all 
algebraic properties of their classical specializations.  Their extra 
structure makes them useful for studying their specializations, as it can 
be easier to find patterns before variables are specialized.

The main geometric setting to which these polynomials apply is the 
following.  We have maps of vector bundles 
$$
F_1 \rightarrow  F_2 \rightarrow\dots  \rightarrow F_n  \rightarrow 
E_n \rightarrow  \dots  \rightarrow E_2 \rightarrow E_1 \tag 1
$$
on a variety or scheme  $X$, where each $F_i$ and $E_i$ has rank $i$.  
We do {\it not\/} assume here that the maps $F_i \rightarrow F_{i+1}$ 
are injective, or that the maps $E_{i+1} \rightarrow E_i$ are 
surjective, as was the case studied in \cite{F2}.  For each $w$ in the 
symmetric group $S_{n+1}$ there is a degeneracy locus
$$
\Omega_w = \{ x \in X \mid \operatorname{rank}\, (F_q(x) \rightarrow 
E_p(x)) \leq r_w(p,q) \text{ for all }  1 \leq p,q \leq n\}, \tag 2
$$
where $r_w(p,q)$ is the number of $i \leq p$ such that $w(i) \leq q$.  
Such degeneracy loci will be described by the double form 
$\frak{S}_w(c,d)$ of universal Schubert polynomials, evaluated at the 
Chern classes of all the bundles involved.  Unlike the situation studied 
in \cite{F2}, where these Chern classes were determined by their first 
Chern classes, in the present general setting one must have more 
general polynomials to describe such loci.  There are similar formulas 
when some of the bundles in (1) are missing.

All Schubert polynomials are indexed by permutations $w$ in some 
symmetric group $S_{n+1}$.  We will present these universal Schubert 
polynomials in two forms.  The first (in its single form), denoted 
$\frak{S}_w(c)$, is a polynomial in variables  $c_i(k)$, for $1 \leq i 
\leq k \leq n$.  When $c_i(k)$ is specialized to the $i^{\text{th}}$ 
elementary symmetric polynomial $e_i(x_1, \dots, x_k)$ in variables 
$x_1, \dots ,x_k$, this polynomial becomes the classical Schubert 
polynomial, denoted $\frak{S}_w(x)$.  When $c_i(k)$ is specialized to 
the  $i^{\text{th}}$ quantum elementary symmetric polynomial, in
variables 
$x_1, \dots, x_k, q_1, \dots, q_{k-1}$, $\frak{S}_w(c)$ specializes to 
the quantum Schubert polynomial  $\frak{S}_w^q$ of \cite{F-G-P}.  In 
our geometric setting, $c_i(k)$  will be the $i^{\text{th}}$ Chern class 
of a vector bundle of rank $k$.

We write the second form of the universal Schubert polynomial, 
denoted $\frak{S}_w(g)$,  as a polynomial in variables $g_i[j]$, for 
$i \geq 1$ and $j \geq 0$ with $i+j \leq n+1$; we regard $g_i[j]$ as an 
indeterminate of degree $j+1$.  This polynomial $\frak{S}_w(g)$ is 
obtained from $\frak{S}_w(c)$ by replacing each  $c_i(k)$ by the 
coefficient of $T^i$ in the determinant of $A+IT$, where $A$ is the $k$ 
by $k$ matrix with $g_i[j \! - \! i]$ in the $(i,j)$ position for $i \leq j$, 
and with $-1$ in positions $(i+1,i)$ below the diagonal, and $0$ elsewhere.  
(See Section 4 
for another definition of these polynomials.)  One recovers the classical 
Schubert polynomials $\frak{S}_w(x)$ by setting $g_i[0] = x_i$ and 
$g_i[j] = 0$ for $j \geq 1$, and one recovers the quantum Schubert 
polynomials $\frak{S}_w^q$ by setting $g_i[0] = x_i$,  $g_i[1] = q_i$, 
and $g_i[j] = 0$ for $j \geq 2$.  In Section 4 we will see that other 
specializations give the polynomials defined in \cite{CF2} for Schubert 
classes in quantum cohomology rings of partial flag varieties.  Since 
the variables $c_i(k)$ and $g_i[j]$ generate the same polynomials ring, 
i.e., $\Bbb{Z}[c] = \Bbb{Z}[g]$, the two forms of universal polynomials 
are equivalent.

From the single polynomials we will construct universal double 
Schubert polynomials $\frak{S}_w(c,d)$ and $\frak{S}_w(g,h)$, which 
also specialize to the known cases of double Schubert polynomials.  
In the geometric setting of \thetag{1}, 
the variables $c_i(j)$ become the Chern classes $c_i(E_j)$, and the 
variables $d_i(j)$ become $c_i(F_j)$.

The universal polynomials are constructed in Section 2.  The theorems 
relating them to degeneracy loci are proved in Section 3.  The 
last section contains some determinantal formulas for universal 
Schubert polynomials, and some results and questions about their 
algebra.

In \cite{F2}, following the classical approaches of \cite{B-G-G} and 
\cite{D}, the degeneracy loci formulas were proved -- in a universal 
setting on a flag bundle -- by starting with the locus of top 
codimension, which is realized as the zero of a section of a vector 
bundle; then the other loci are constructed inductively by a sequence of 
$\Bbb{P}^1$--bundle correspondences.  In the present setting, neither 
of these methods is available.  Indeed, the top double classes 
$\frak{S}_w(c,d)$ do not factor.  Our procedure, roughly speaking, is to 
find a locus in a flag bundle that maps to a given degeneracy locus 
$\Omega_w$, but where one has injections and surjections of the 
bundles, so that one can apply the results of \cite{F2}; then this 
formula is pushed forward to get a formula for  $\Omega_w$.  

Ionu\c t  Ciocan-Fontanine initiated this project by asking several 
years ago for a degeneracy locus formula that would apply when the 
maps $E_{i+1} \to E_i$ are not surjective.  His work in \cite{CF1}, 
\cite{CF2} was a source for this question, and conversations with him 
have been very useful.  Algebraically, the universal Schubert 
polynomials are natural generalizations of the quantum Schubert 
polynomials of \cite{F-G-P}, cf\. \cite{K-M} and \cite{CF-F}, and the 
inspiration of these sources should be clear.  I thank Chandler 
Fulton and Mel Hochster for advice for computer testing related to this
study.

The fact that the universal versions of these degeneracy loci are Cohen-
Macaulay was conjectured in the preprint version of this paper.  Claudio
Procesi pointed out the paper \cite{A-DF-K}, which proved, in
characteristic zero, that the reduced structures on these loci are Cohen-
Macaulay.  Recently Lakshmibai and Magyar \cite{L-M} succeeded in
proving that these schemes, with their natural determinantal structures,
are reduced and Cohen-Macaulay, in all characteristics.  This strengthens
our results, so that the formulas for the degeneracy loci have the usual
meaning in intersection theory, as in \cite{F1, \S 14}.

\head 2. Definitions of Universal Schubert Polynomials \endhead

We will give three constructions of the single universal Schubert 
polynomials $\frak{S}_w(c)$.  We consider independent variables 
$c_i(j)$ for $1 \leq i \leq k \leq n$.  It is to be understood that 
$c_i(j) = 1$ if $i = 0$, and $c_i(j) = 0$  if $i < 0$ or $i > j$.  Let $w$ be 
a permutation in $S_{n+1}$, and let $l(w)$ denote its length.
 
The quickest definition of $\frak{S}_w(c)$ is a variation of that used in 
\cite{F-G-P} to define quantum Schubert polynomials. In the classical
case such a formula can be found in \cite{L-S2, (2.10)}.  A classical 
Schubert polynomial $\frak{S}_w(x)$ can be written uniquely in the 
form 
$$
\frak{S}_w(x) = \sum a_{i_1, \dots, i_n}e_{i_1}(x_1)\cdot \dots \cdot 
e_{i_n}(x_1, \dots, x_n),
$$
the sum over  $(i_1, \dots, i_n)$ with each $i_\alpha \leq \alpha$ and 
$\sum i_\alpha = l(w)$; here $a_{i_1, \dots, i_n}$ are unique integers 
(depending on $w$).  Define
$$
\frak{S}_w(c) = \sum a_{i_1, \dots, i_n} \, c_{i_1}(1)\cdot \dots \cdot 
c_{i_n}(n). \tag 3
$$

The preceding definition is based on the following elementary fact (see 
\cite{L-S2, (2.7)} and \cite{F-G-P, Prop. 3.4}), that we will 
use frequently.

\proclaim{Lemma 2.1} 
Let $R$ be a commutative ring, and let  $M$ be the free $R$-submodule 
of the polynomial ring $R[c]$ spanned by all monomials 
$c_{i_1}(1)\cdot \dots \cdot c_{i_n}(n)$ with each $i_\alpha \leq 
\alpha$.  Let $M'$ be the free $R$-submodule of the polynomial ring 
$R[x]$ spanned by all monomials $x_1^{j_1}\cdot \dots \cdot x_n^{j_n}$ 
with each $j_\alpha \leq n+1-\alpha$.  Then the map which sends 
$c_i(j)$ to $e_i(x_1, \dots ,x_j)$ determines an isomorphism of $M$ 
onto $M'$.
\endproclaim

Universal Schubert polynomials can also be defined by a direct 
inductive procedure analogous to that for the classical Schubert 
polynomials.   For  $w_0$ the permutation in $S_{n+1}$ of longest 
length, i.e., $w_0(i) = n+2-i$ for $1 \leq i \leq n+1$, set  
$$
\frak{S}_{w_0}(c) = c_1(1) \cdot c_2(2) \cdot \dots \cdot c_n(n).
\tag 4
$$
The general universal polynomial is determined by the property that if 
$k$ is an integer with $w(k) < w(k+1)$, and $v = w @,@, s_k$ is obtained 
from $w$ by interchanging the values of $k$ and $k+1$, then 
$$
\frak{S}_w(c) = \partial_k(\frak{S}_v(c)). \tag 5
$$ 
Here $\partial_k$ is an additive endomorphism of the free $\Bbb{Z}$-
module $M$ spanned by monomials $c_{i_1}(1)\cdot \dots \cdot 
c_{i_n}(n)$ with each $i_\alpha \leq \alpha$.  It takes such a monomial 
to a sum of signed monomials, each with the same indices except in 
positions $k-1$ and $k$.  Write for simplicity $[a,b]$ for the monomial 
with $c_a(k-1)$ and  $c_b(k)$ in these two positions, with the other 
positions fixed but arbitrary.  Note that $[p,q] = 0$ if $p$ or $q$ is 
negative, or if $p > k-1$ or $q > k$.  With this notation the formula for 
$\partial_k$ is
$$ \alignedat2
\partial_k([a,b]) &= \sum_{i \geq 0} \, [a+i,b-1-i] - 
\sum_{i \geq 1} \, [b-1-i,a+i]  &&\qquad  \text {if } \, \, a \geq b-1; \\
 \partial_k([a,b]) &= \sum_{i \geq 0} \, [b-1+i,a-i] - 
\sum_{i \geq 1} \, [a-i,b-1+i]  &&\qquad  \text {if } \, \, a \leq b-2.
\endalignedat \tag 6
$$
For $k = 1$, $\partial_1([1]) = [0]$ and $\partial_1([0]) = 0$. See
\cite{L-S2, \S3} for similar formulas.

To see that this definition is well-defined, it suffices to verify that 
$\partial_k \circ \partial_l = \partial_l \circ \partial_k$ if $| k-l | 
\geq 2$, and that $\partial_k \circ \partial_{k+1} \circ \partial_k = 
\partial_{k+1} \circ \partial_k \circ \partial_{k+1}$ for $1 \leq k \leq 
n-1$.  This can be verified directly from the definition, but it follows 
easily from the lemma, together with the simple verification that 
$\partial_k$ is compatible with the standard difference operator 
$\partial_k^{(x)}$ defined on $R[x]$: 
$\partial_k^{(x)}(P) = (P - s_k(P))/(x_k - x_{k+1})$,  where  $s_k(P)$  is 
the result of interchanging $x_k$ and $x_{k+1}$ in $P$ (cf\. \cite{M], 
Chapter 2}).  It also follows from this argument that the universal 
quantum Schubert polynomials form a basis for $M$, but we will see a 
stronger reason for this in Proposition 2.2.

The third definition defines double versions $\frak{S}_w(c,y)$ of these 
polynomials, with general first variables $c_i(j)$ as above and special 
second variables $y_1, \dots, y_n$.  This definition is similar to a 
construction in \cite{K-M}, cf\. \cite{CF-F}, except that we have 
included signs with the $y$ variables in order to be consistent with the 
original notation of Lascoux and Sch\"utzenberger.  For this, set
$$
\frak{S}_{w_0}(c,y) = \prod_{i = 1}^n \bigl(\,\sum_{j=0}^i c_{i-j}(i)
(-y_{n+1-i})^j \, 
 \bigr). \tag 7
$$
Now if  $k+1$ appears in the list of values of $w$ to the right of $k$, 
and  $v$ interchanges the positions of $k+1$ and $k$, i.e., $v = s_k @,@, 
w$,  then
$$
\frak{S}_w(c,y) = - \, \partial_k^{(y)}(\frak{S}_v(c,y)), \tag 8
$$ 
where $\partial_k^{(y)}$ is the standard difference operator, acting on 
the $y$ variables alone.  Then 
$$
\frak{S}_w(c) = \frak{S}_w(c,0), \tag 9
$$
i.e., the single universal Schubert polynomials are obtained from these 
double polynomials by setting the second set of variables equal to $0$.

It is not hard to see that these three definitions agree.  That the first 
and second definitions agree follows from the fact that they both give 
polynomials in the module $M$ of the lemma which specialize to the 
classical Schubert polynomials under the isomorphism from $M$ to 
$M'$.  Next we observe that the double polynomials $\frak{S}_w(c,y)$ 
specialize to the usual double Schubert polynomials $\frak{S}_w(x,y)$ 
when each $c_i(j)$ is sent to $e_i(j)$.  This follows from the fact that 
$\frak{S}_w(x,y) = (-1)^{l(w)}\frak{S}_{w^{-1}}(y,x)$, \cite{M, 
(6.4)(iii)}.  From this it follows that  $\frak{S}_w(c,0)$ specializes to 
$\frak{S}_w(x) = \frak{S}_w(x,0)$ when $c_i(j)$ is replaced by 
$e_i(j)$, and this shows that the third definition agrees with the first 
two. 

There is a natural definition of {\it universal double Schubert 
polynomials}, that we denote by  $\frak{S}_w(c,d)$,  where  $c$ stands 
for the variables $c_i(j)$ and $d$ stands for another set of variables 
$d_i(j)$.  These are defined by the formula
$$
\frak{S}_w(c,d) = \sum_{u,v} (-1)^{l(v)}\frak{S}_u(c)\frak{S}_v(d), \tag 10
$$
where the sum is over all  $u$ and $v$ in $S_{n+1}$ such that 
$v^{-1}\, u = w$ and $l(u)+l(v) = l(w)$.  The same argument as in the 
preceding paragraph, together with \cite{M, (6.3)}, shows that 
$\frak{S}_w(c,d)$ specializes to the polynomials $\frak{S}_w(c,y)$ 
when each $d_i(j)$ is specialized to the $i^{\text{th}}$ elementary 
symmetric polynomial in $y_1, \dots, y_j$.  In particular, when 
$c_i(j)$ is also specialized to $e_i(x_1, \dots,x_j)$, then  
$\frak{S}_w(c,d)$ becomes the classical double Schubert polynomial 
$\frak{S}_w(x,y)$ of \cite{L}, cf. \cite{M}.

The universal double Schubert polynomials for permutations in $S_3$ are:
$$
\align
\frak{S}_{3@,@,2@,@,1} &= c_1(1)c_2(2)-(c_1(1)c_1(2)-c_2(2))d_1(1)-
c_2(2)d_1(2) \\
& \qquad +c_1(1)(d_1(1)d_1(2)-d_2(2))+c_1(2)d_2(2)-d_1(1)d_2(2) \\
\frak{S}_{2@,@,3@,@,1} &= c_2(2)-c_1(2)d_1(1)+d_1(1)d_1(2)-d_2(2) \\
\frak{S}_{3@,@,1@,@,2} &= c_1(1)c_1(2)-c_2(2)-c_1(1)d_1(2)+d_2(2) \\
\frak{S}_{1@,@,3@,@,2} &= c_1(2)-d_1(2) \\
\frak{S}_{2@,@,1@,@,3} &= c_1(1)-d_1(1) \\
\frak{S}_{1@,@,2@,@,3} &= 1
\endalign
$$
One can make specializations of either or both variables.  For example, 
introducing variables $h_i[j]$ analogous to the variables $g_i[j]$, one 
has polynomials that we denote by $\frak{S}_w(g,h)$.  In general we let 
the position in the alphabet distinguish among these different single 
and double universal Schubert polynomials, using variables $c,\, d$ for 
the first kind, variables $g,\, h$ for the second, and $x,\, y$ for the 
classical case.  (This seems preferrable to introducing different 
notations for each realization of these polynomials.)

By construction the universal Schubert polynomials are expressed as a 
linear combination of monomials $c_{i_1}(1)\cdot \dots \cdot 
c_{i_n}(n)$.  We want to say a little more about this expansion.  For 
this we need a modified version of the code of a permutation $w$.  For 
$w \in S_{n+1}$, we will define this code, and denote it $c'(w)$, to be 
the sequence $(i_1, \dots, i_n)$, where $i_k$ is defined by the formula
$$
i_k = \operatorname{Card} \{\, j \leq k \mid w(j) > w(k+1)\,\}. \tag 11
$$
This number $i_k$ is the number of boxes in the $k^{\text{th}}$ row of a 
modified diagram $D'(w)$ of $w$, which is constructed as follows.  
Form an $n$ by $n$ square of boxes arranged as in a matrix, and, for 
$2 \leq i \leq n+1$, remove all the boxes from the row directly above 
the position $(i,w(i))$ that are strictly to the left of this position, and 
also all the boxes in the column directly to the left of and strictly 
above that position. So $D'(w)$ consists of boxes $(i,j)$ such that 
$w(i+1) \leq j$ and $w^{-1}(j+1) \leq i$.  It is easily seen that 
if $c'(w) = (i_1, \dots,i_n)$, 
then the Lehmer code of the permutation $w_0\, w \, w_0$ (cf\. 
\cite{M, p. 9}) is $(i_n, \dots, i_1, 0)$.  In particular, this code $c'(w)$ 
determines $w$, and the sum of the integers in $c'(w)$ is the length of 
$w$ -- properties that are easily proved directly.  
 
\proclaim{Propositon 2.2} 
Let $c'(w) = (i_1, \dots,i_n)$.  Then
$$
\frak{S}_w(c) = c_{i_1}(1)\cdot \dots \cdot c_{i_n}(n) +
\sum n_{j_1, \dots,j_n} \, c_{j_1}(1)\cdot \dots \cdot c_{j_n}(n),
$$
where the sum is over $(j_1, \dots,j_n)$ that are strictly smaller 
than  $(i_1, \dots,i_n)$ in the lexicographic ordering.
\endproclaim

\demo{Proof} This will follow from the second construction of the 
universal Schubert polynomials.  The assertion is trivial when $w = 
w_0$, so we may assume it for all $v$ of length greater than the length 
of a given $w$.  Let $k$ be the smallest integer such that $w(k) < 
w(k+1)$, and let $v = w@,@, s_k$.  Since $\frak{S}_w(c) = 
\partial_k(\frak{S}_v(c))$, it suffices to show that the assertion of 
the proposition for $v$ implies the assertion for $w$.  When $k = 1$, 
this is completely straightforward, since $\partial_k$ is so simple in 
this case.  
For $k > 1$, the code $I = c'(w)$ has $i_j = j$ for $1 \leq j \leq k-1$.  
The code $H = c'(v)$ is that same as $I$ except in positions $k-1$ and 
$k$, where it is $i_k$ and $k$ respectively.  Write $c_J$ for 
$c_{j_1}(1)\cdot \dots \cdot c_{j_n}(n)$.  It suffices to check that 
$\partial_k(c_H) = c_I \pm \text{smaller terms}$, and that, if $J < H$, 
then $\partial_k(c_J)$ consists entirely of terms that are smaller that 
$c_I$.  The verification of these facts is straightforward from 
\thetag{6}.
\enddemo

Note that this proposition implies (and is equivalent to) the 
corresponding assertion for the expression of the classical Schubert 
polynomials in terms of elementary symmetric monomials, or the 
quantum Schubert polynomials in terms of quantum elementary 
symmetric polynomials.  It shows effectively why the universal 
Schubert polynomials form a basis for the module $M$ that appears in 
Lemma 2.1.

Another property, which follows immediately from the definitions, is 
the {\it stability} of these universal polynomials: if $i$ is the 
canonical embedding of $S_{n+1}$ in $S_{n+2}$, then 
$\frak{S}_{i(w)}(c,d) = \frak{S}_w(c,d)$.  Thus single and double 
universal polynomials are defined for $w$ in $S_{\infty} = \cup S_n$.

Finally, we have the expected duality property, which is also an 
immediate consequence of the definitions:
$$
\frak{S}_w(d,c) = (-1)^{l(w)}\frak{S}_{w^{-1}}(c,d). \tag 12
$$

\head  3. Formulas for Degeneracy Loci \endhead

In this section we explain how the universal double Schubert 
polynomials describe degeneracy loci of appropriate maps of vector 
bundles.  We assume that we are given the situation of vector bundles 
and maps as described in \thetag{1} of the introduction, on an algebraic 
scheme $X$ over a field.  (This is only to simplify the exposition; the 
procedures of \cite{F1, \S 20} show how to modify the arguments 
for schemes of finite type over an arbitrary regular base scheme.)  Set 
$c_i(j) = c_i(E_j)$, the $i^{\text{th}}$ Chern class of $E_j$, and set 
$d_i(j) = c_i(F_j)$.  We claim that $\frak{S}_w(c,d)$ is the formula for 
the locus $\Omega_w$ defined in \thetag{2}.  

As usual, of course, this assertion must be interpreted correctly, 
depending on assumptions about how general the maps between the 
vector bundles are.  When $X$ is an algebraic scheme of pure dimension
$k$, there is a class $\boldsymbol{\Omega}_w$ in the Chow group 
$A_{k-l(w)}(\Omega_w)$, whose image in $A_{k-l(w)}(X)$ is
$\frak{S}_w(c,d) {\ssize{\cap}} [X]$.  When $\Omega_w$ has the expected
dimension $k-l(w)$, then $\boldsymbol{\Omega}_w$ is a positive cycle
whose support is $\Omega_w$.  If, in addition, $X$ is Cohen-Macaulay,
then $\boldsymbol{\Omega}_w = [\Omega_w]$, where $\Omega_w$ has its
natural subscheme structure defined locally by vanishing of determinants. 
More generally, $\boldsymbol{\Omega}_w = [\Omega_w]$ whenever 
$\operatorname{depth}(\Omega_w,X) = l(w)$.  Without any depth
conditions, $[\Omega_w]$ can have larger multiplicities than
$\boldsymbol{\Omega}_w$.  The construction of the class
$\boldsymbol{\Omega}_w$ is compatible with proper pushforward, and
pullback by flat or l.c.i. morphisms.  In fact, the class
$\boldsymbol{\Omega}_w$ can be constructed to live in the bivariant
group $A^{l(w)}(\Omega_w \to X)$; for complex varieties, similar
constructions produce classes in the relative cohomology groups 
$H^{2l(w)}(X,\,X \smallsetminus \Omega_w)$.  This whole package is what
is meant by the phrase ``$\frak{S}_w(c,d)$ gives a formula for
$\Omega_w$.''  For details about these interpretations, see 
\cite{F1, \S 14 and \S 17}, and \cite{F-P, App. A}.

The construction of $\boldsymbol{\Omega}_w$ is carried out as follows. 
Let
$$
H = \bigoplus_{i=1}^{n-1}\operatorname{Hom}(F_i,F_{i+1}) \oplus
\operatorname{Hom}(F_n,E_n) \oplus 
\bigoplus_{i=1}^{n-1}\operatorname{Hom}(E_{i+1},E_i),
$$
regarded as a bundle over  $X$.  On $H$ there are universal or tautological
maps between the pullbacks of the bundles, and hence there are universal
loci $\widetilde{\Omega}_w \subset H$.   If $X$ is reduced, it follows
from the theorem of Lakshmibai and Magyar \cite{L-M} that these schemes
$\widetilde{\Omega}_w$ are reduced, of codimension $l(w)$.  Moreover, if
$X$ is Cohen-Macaulay, their theorem implies that 
$\widetilde{\Omega}_w$ is Cohen-Macaulay.  Note that on any open set
$U$ of $X$ on which the bundles are trivial, $H$ is a product $U \times
H_0$, where $H_0$ is the corresponding bundle constructed from vector
spaces, and $\widetilde{\Omega}_w$ is the product of $U$ and the
corresponding universal locus in $H_0$; it is these universal local loci
that are studied in \cite{L-M}.  

The given maps on $X$ determine a section $s \! : \! X \to H$, and
$\Omega_w = s^{-1}(\widetilde{\Omega}_w)$.  We define
$\boldsymbol{\Omega}_w$ to be the refined pullback of the class of
$\widetilde{\Omega}_w$, in the sense of \cite{F1}, i.e.,
$$
\boldsymbol{\Omega}_w = s^{\,!\,}[\widetilde{\Omega}_w].
$$
To prove that $\frak{S}_w(c,d)$ is the formula for 
the locus $\Omega_w$, it suffices to prove that
$$
[\widetilde{\Omega}_w] = \frak{S}_w(c,d)\, {\ssize{\cap}} /, [H]
\tag 13 
$$
in $A_{\dim(H)-l(w)}(H)$.  As in \cite{F1, \S 14} and \cite{F-P, App. A},
the fact that the univeral local loci are Cohen-Macaulay is what makes this
construction have all the stated properties.  In the proofs that follow, we
therefore replace $X$ by $H$, so we may assume the maps are locally
universal in this sense. 

When the maps $F_i \to F_{i+1}$ are injective, and the maps $E_{i+1} \to
E_i$ are surjective, the formula is exactly that proved in \cite{F2}. 

Now consider the situation where we do not assume the maps from 
each $F_i$ to $F_{i+1}$ are injective, but for now we assume the maps 
$E_{i+1} \to E_i$ are surjective.  

\proclaim{Proposition 3.1} The formula for $\Omega_w$  is
$\frak{S}_w(x,d)$, where $x_i =   
c_1(\operatorname{Ker}(E_i\to E_{i-1})$ and $d_i(j) = c_i(F_j)$.
\endproclaim
 \demo{Proof}
We may assume that the map from $F_n$ to $E_n$ factors into an
inclusion $F_n \to V$ followed by a surjection $V \to E_n$, where $V$ is a
vector bundle of rank $n+1$.  Indeed, as in \cite{F2}, one first considers
the factorization $F_n \to F_n \oplus E_n \to E_n$ given by the graph of
$F_n \to E_n$, and then pulls back to flag bundles to fill in between $F_n$
and $V = F_n \oplus E_n$, and between $V$ and $E_n$; this replaces $n+1$
by $2n$.  

Let $\rho  \! : \!  F \to X$ be the bundle of complete flags in $V$, with
universal subbundles $U_i$, so we have 
$$
U_1 \subset  \dots \subset U_n \subset \rho^*(V) \rightarrow 
\rho^*(E_n) \rightarrow \dots \rightarrow \rho^*(E_1).
$$
From these bundles and maps, by the case just considered, for each $w$ 
in $S_{n+1}$, we have a locus $\widetilde\Omega_w$ on $F$, and we know
that
$$
[\widetilde\Omega_w] = \frak{S}_w(x,y) \, {\ssize{\cap}}/, [F],
$$
where $y_i = c_1(U_i/U_{i-1})$.

Let $Z_n \subset F$ be the locus given by the vanishing of the canonical 
map from $\rho^*(F_n)$ to $V/U_n$.  This is given by the vanishing of a 
section of the bundle $\rho^*(F_n)^\vee \otimes V/U_n$, so 
$$
[Z_n] = c_n(\rho^*(F_n)^\vee \otimes V/U_n) \, {\ssize{\cap}}/, [F].
$$
Notice that on 
$Z_n$ the map from $\rho^*(F_n)$ to $V$ factors through $U_n$.  On 
$Z_n$ we have the locus $Z_{n-1}$ given by the vanishing of the map 
from $\rho^*(F_{n-1})$ to $U_n/U_{n-1}$, so $Z_{n-1}$ is represented on
$Z_n$ by a top Chern class.  On $Z_{n-1}$ the map from 
$\rho^*(F_{n-1})$ to $U_n$ factors through $U_{n-1}$.  Continuing in 
this way we get a sequence of loci $Z_1 \subset \dots  \subset Z_n 
\subset F$.  Let $Z = 
Z_1$.  Then
$$
[Z] = \prod_{i=1}^n \bigl( c_i(\rho^*(F_i)^\vee \otimes U_{i+1}/U_i \bigr)
\, {\ssize{\cap}}\, [F].
$$
It follows from the fact that the given maps are locally universal that
$\widetilde\Omega_w$ meets $Z$ properly, and that $\rho$ maps 
$\widetilde\Omega_w \cap Z$ birationally onto $\Omega_w$.  (See the
remark following the proof.)  It therefore suffices to show that
$$
\rho_* \bigl( \, \prod_{i=1}^n ( c_i(\rho^*(F_i)^\vee \otimes U_{i+1}/U_i) 
\cdot  \frak{S}_w(x,y) \bigr) = \frak{S}_w(x,d), \tag 14
$$
where $\rho_*$ is the pushforward from $A^{m+l(w)}(F)$ to 
$A^{l(w)}(X)$, with $m = n(n+1)/2$.  Let $Y_k$ be the flag 
bundle of subbundles of $V$ of all ranks from $k$ through $n$.  So 
$\rho$ factors into a composite of projective bundle projections
$$
F = Y_1 \to Y_2 \to \dots \to Y_n = \Bbb{P}^*(V) \to Y_{n+1} = X.
$$
On $Y_k$ we have universal bundles $U_k \subset \dots \subset U_n 
\subset V$.  (Here and elsewhere in the proof we use the same notation 
for bundles as for their pullbacks by canonical projections.)  We know 
that we can write 
$$
\frak{S}_w(x,y) = \sum a_I(x)c_{i_1}(U_1)c_{i_2}(U_2) \dots 
c_{i_n}(U_n),
$$
where the sum is over $I = (i_1, \dots, i_n)$ with $i_\alpha \leq 
\alpha$, and the $a_I(x)$ are polynomials in the $x$ variables.  (Note 
for this that $c_j(U_i) = e_j(y_1, \dots , y_i)$.)  It suffices to prove 
that for each such $I$, and for each $2 \leq k \leq n+1$, the 
pushforward of
$$
\prod_{i=1}^n c_i(F_i^\vee \otimes U_{i+1}/U_i) \cdot c_{i_1}(U_1) 
\dots c_{i_n}(U_n)
$$
from $F = Y_1$  to $Y_k$ is equal to
$$
\prod_{i=k}^n c_i(F_i^\vee \otimes U_{i+1}/U_i) \cdot c_{i_1}(F_1) 
\dots c_{i_{k-1}}(F_{k-1}) \cdot c_{i_k}(U_k) \dots c_{i_n}(U_n).
$$  
To prove this, using the projection formula, it suffices by induction to 
show that the pushforward from $Y_k$ to $Y_{k+1}$ of 
$c_k(F_k^\vee \otimes U_{k+1}/U_k) \cdot c_i(U_k)$ is equal to 
$c_i(F_k)$.  But this is a special case of the following elementary 
Gysin formula for projective bundles.
\enddemo

\proclaim{Lemma 3.2} Let $G$ be a vector bundle of rank $k+1$ on a 
scheme $Y$, let $\Bbb{P} = \Bbb{P}^*(G)$ be the projective bundle of 
hyperplanes in $G$, with projection $p \! : \! \Bbb{P} \to Y$ and 
tautological sequence
$$
0 \to H \to p^*(G) \to \Cal{O}(1) \to 0
$$
of bundles on $\Bbb{P}$.  Let $K$ be a vector bundle of rank $k$ on $Y$.  
Then, for $0 \leq i \leq k$, 
$$
p_* \bigl( c_k(p^*(K)^\vee \otimes \Cal{O}(1)) \cdot c_i(H)  \bigr) = 
c_i(K).
$$
\endproclaim
\demo{Proof}
Let $\zeta = c_1(\Cal{O}(1))$.  We use the basic fact that for any 
integer $r$,  $p_*(\zeta^{r+k}) = (-1)^rs_r(G)$, where $s_r(G)$ denotes 
the $j^{\text{th}}$ Segre (inverse Chern) class of the bundle $G$.  We 
use the formula for the top Chern class of a tensor product with a line 
bundle
$$ 
c_k(p^*(K)^\vee \otimes \Cal{O}(1)) = 
\sum_{a=0}^k (-1)^a c_a(p^*K)\zeta^{k-a},
$$
and the Whitney formula 
$c_i(H) = \sum_{b=0}^i (-1)^bc_{i-b}(p^*G)\zeta^b$.  
This gives
$$
p_* \bigl( c_k(p^*(K)^\vee \otimes \Cal{O}(1)) \cdot c_i(H)  \bigr) = 
\sum_{a=0}^k (-1)^a c_a(K) \bigl( \sum_{b=0}^i (-1)^{b+b-a} 
c_{i-b}(G)s_{b-a}(G) \bigr).
$$
The inner sum vanish unless $a = i$, when it gives $(-1)^i$, so the right 
side is the required $c_i(K)$.
\enddemo

This completes the proof of Proposition 3.2, which constructs the 
required class in the case where the right maps are surjective but the 
left maps are not necessarily injective.  The dual case, where the left 
maps are injective and the right maps are arbitrary, can be handled by a 
dual construction, or simply by taking the duals of all the maps, and 
applying the case just considered to the situation
$$
E_1^\vee \rightarrow \dots \rightarrow  E_n^\vee \rightarrow F_n^\vee 
\rightarrow \dots \rightarrow  F_1^\vee.
$$
The locus $\Omega_w$ for the original maps is the locus 
$\Omega_{w^{-1}}$ for this dual sequence.  So we define the class 
$\boldsymbol{\Omega}_w$ for the original maps to be the class 
$\boldsymbol{\Omega}_{w^{-1}}$ for this dual.  The duality property 
\thetag{12} guarantees that this class has the right image in the Chow 
group of $X$.

\example{Remark 3.3} In the preceding proof we used the fact that the
scheme $\widetilde{\Omega}_w \cap Z$ is reduced when $X$ is reduced.  Although this
can be seen directly, this argument can be avoided.  For if not, the
argument would prove that $\frak{S}_w(x,d)$ represents some multiple of
the class of the reduced subscheme of $\Omega_w$.  This would imply
that the polynomial $\frak{S}_w(x,d)$ is a power of some other
polynomial, but we know that this is not the case even after specializing
to $\frak{S}_w(x,0)$.
\endexample

\proclaim{Proposition 3.4} Suppose we are given, on a scheme $X$, 
vector bundles and maps 
$$
F_1 \rightarrow  \dots  \rightarrow F_n  \rightarrow E_n \rightarrow  
\dots  \rightarrow E_1,
$$
with the ranks of $E_i$ and $F_i$ being $i$.  Then for each $w$ in 
$S_{n+1}$ the formula for $\Omega_w$ is the 
universal double Schubert polynomial $\frak{S}_w(c,d)$, where $c_i(j) 
= c_i(E_j)$ and $d_i(j) = c_i(F_j)$.
\endproclaim
\demo{Proof}
We proceed exactly as in the preceding proposition, first reducing to 
the case where one has an intermediate bundle $V$.  One then has the 
formula for the corresponding locus $\widetilde{\Omega}_w$ on the flag
bundle $F$, this time applied to the case where the second maps are
arbitrary.  This class maps to $\frak{S}_w(c,y) {\ssize{\cap}} [F]$ in the
Chow group of $F$, where 
$c_i(j) = c_i(E_j)$ and $y_i = c_1(U_i/U_{i-1})$.  We apply the same 
construction as before, pushing down the product of this class and the 
class $[Z]$.  This time we must prove the formula
$$
\rho_* \bigl( \, \prod_{i=1}^n ( c_i(\rho^*(F_i)^\vee \otimes U_{i+1}/U_i) 
\cdot  \frak{S}_w(c,y) \bigr) = \frak{S}_w(c,d).
$$
The only difference is that the polynomials $a_I(x)$ that 
appeared in the preceding proof are replaced by polynomials in the classes 
$c_i(E_j)$ for this proof.
\enddemo

\example{Remark 3.5} As the proof shows, the Schubert polynomials 
$\frak{S}_w(c,y)$ represent the loci $\Omega_w$ for the situation of 
the theorem, but under the conditions that the maps $F_i \to F_{i+1}$ 
are injective.  Here $c_i(j) = c_i(E_j)$, and $y_i = c_1(F_i/F_{i-1})$.
One has a dual interpretation for $\frak{S}_w(x,d)$. \endexample

There is a corresponding and more general theorem for maps of bundles 
of arbitrary increasing and decreasing ranks.  Suppose we have vector 
bundles and maps
$$
F_1 \rightarrow  \dots  \rightarrow F_s  \rightarrow E_r \rightarrow
\dots  \rightarrow E_1,
$$
with $\operatorname{rank}(E_i) = a_i$,  $\operatorname{rank}(F_i) = 
b_i$, such that 
$$
a_1 < a_1 < \dots < a_r  \qquad \text {and} \qquad b_1 < b_1 < \dots < 
b_s.
$$  
Let $w$ be a permutation whose diagram $D'(w)$ is contained in
$A \times B$, with $A = \{a_1, \dots , a_r \}$ and $B =
\{b_1, \dots , b_s \}$.  There  
is a degeneracy locus $\Omega_w \subset X$ 
defined by the conditions that the rank of the map from $F_j$ to $E_i$ 
is at most $r_w(a_i,b_j)$ for all $i$ and $j$. 

\proclaim{Lemma 3.6} If $D'w) \subset A \times B$, then no $c_i(k)$
occurs in $\frak{S}_w(c,d)$ with $k \notin A$, and no $d_j(l)$ occurs in
$\frak{S}_w(c,d)$ with $l \notin B$. \endproclaim
\demo{Proof} This follows from two claims, valid for any permutation
$w$ with
code $c'(w) = (i_1, \dots , i_n)$ and an integer $k$  such that 
$i_k = 0$. The first claim is that all terms $c_{j_1}(1)\cdots c_{j_n}(n)$
occurring in $\frak{S}_w(c)$ have $j_k = 0$; this is proved as in the
proof of Proposition 2.2.  The second claim is that if $i$ is an integer
with $w^{-1}(i+1) < w^{-1}(i)$, and $v = s_i @,@, w$, and $c'(w) =
(j_1, \dots,j_n)$, then $j_k = 0$; this follows immediately from the 
definition.  It then follows that in the sum \thetag{10} expressing
$\frak{S}_w(c,d)$, all the $\frak{S}_u(c)$ have no term $c_i(k)$ with
$k \notin A$.  The assertion for the $d_j(l)$ then follows from 
\thetag{12}.
\enddemo

Let 
$\frak{S}_w(c(E_{\sssize\bullet}),c(F_{\sssize\bullet}))$ denote the 
result of specializing $c_i(a_p)$ to $c_i(E_p)$
and $d_i(b_q)$ to $c_i(F_q)$ in $\frak{S}_w(c,d)$.

\proclaim{Theorem 3.7}  The formula for $\Omega_w$ is 
$\frak{S}_w(c(E_{\sssize\bullet}),c(F_{\sssize\bullet}))$.
\endproclaim
\demo{Proof}
Take $n = \max(r,s)$, and add trivial bundles to the given bundles so 
that one has two sets of bundles of all ranks between $1$ and $n$.  
Insert maps
$$
F_k \to F_k \oplus \bold{1} \to \dots \to F_k \oplus \bold{1}^{b_{k+1} - 
b_k -1} \to F_{k+1},
$$
where the maps to successive additions of trivial factors are the 
evident inclusions, and the map from the last bundle to $F_{k+1}$ is the 
given map on the factor $F_k$ and the zero map on the trivial factors.  
The dual construction is made on the other side.  Then one is in a 
position where Proposition 3.4 applies.  The hypotheses on $w$ 
guarantee that all of the rank conditions on the added bundles follow 
from the rank conditions on the given bundles, so that the locus 
$\Omega_w$ is the same whether defined for the given bundles or for 
all the bundles.  Proposition 3.4 gives the required formula.
\enddemo

\example{Remark 3.8} If $w$ is a permutation satisfying the weaker
conditions that the descents of $w$ are contained in $A$ and the 
descents of $w^{-1}$ are contained in $B$, then the same construction as
in the proof of Theorem 3.8 produces a class $\boldsymbol{\Omega}_w$ in
$A_{\dim(X)-l(w)}(\Omega_w)$. Its image in $A_{\dim(X)-l(w)}(X)$ is the
polynomial obtained from $\frak{S}_w(c,d)$ by specializing $c_i(k)$
to $c_i(E_p)$ if $a_p \leq k < a_{p+1}$, with $c_i(k)$ sent to $0$ if
$k < a_1$, and to $c_i(E_r)$ if $k \geq a_r$; similarly $d_i(l)$
is sent to $c_i(F_q)$ if $b_q \leq l < b_{q+1}$.
\endexample

\example{Remark 3.9}  The same class is obtained
if arbitrary bundles are used in place of trivial bundles in the
construction in the proof of the theorem.  However, this is no longer
true for the more general classes constructed in the preceding remark.
\endexample

\example{Remark 3.10} Even in the situation of Proposition 3.4, this 
theorem applies to more general loci $\Omega_w$, where $w$ is not in 
$S_{n+1}$ but $D'(w) \subset [n] \times [n]$.
\endexample

\example{Remark 3.11} The degeneracy loci $\Omega_w$ are defined by 
many rank conditions, but some of them are superfluous.  In fact, the 
rank conditions that are needed on $F_j \to E_i$ are precisely those 
described by the boxes $(a_i,b_j)$ in the diagram we denoted by $D'(w)$ 
in Section 2.  In the case when the first maps were injective and the 
second maps surjective, we saw in \cite{F2} that the essential rank 
conditions were given by those in a subset of the diagram $D(w)$ that 
we called the essential set of $w$.  In fact, this essential set is 
exactly the intersection of $D(w)$ with $D'(w)$.  For general maps as 
considered here, the entire set $D'(w)$ is needed.  The proof of this 
statement is the same as in \cite{F2, Prop. 4.2}.  Note that $D'(w)$ 
consists of $l(w)$ boxes, which is exactly the expected codimension of 
$\Omega_w$.
\endexample

\head 4. Determinantal Formulas \endhead

We first prove a determinantal formula for certain of the polynomials 
$\frak{S}_w(c,y)$, that will be used to show that the universal 
Schubert polynomials specialize to the quantum polynomials for partial 
flag varieties of \cite{CF2}.  For this we fix an integer $l \geq 2$, and 
fix a set
$$
N \, = \, \{ \, n_1 < n_2 < \dots < n_{l-1} < n_l \, \}
$$
of $l$ positive integers.  Set $n_0 = 0$, and set $k_i = n_i - n_{i-1}$ 
for $1 \leq i \leq l$.  Let 
$$
S^{(N)} = \{ \, w \in  S_{n_l} \mid w(i) < w(i+1) \text { if } i \notin N \, 
\}
$$
be the permutations in $ S_{n_l}$ with descents in $N$.  Let $w = 
w_0^{(N)}$ be the element of longest length in $S^{(N)}$, that is, 
$$
w(n_{p-1}+i) = n_l - n_p + i \qquad {\text {for}} \quad 1 \leq i \leq k_p, 
\,\,\, 1 \leq p \leq l.
$$
For nonnegative integers $a, b$, a positive integer $k$, and an arbitrary 
integer $m$, set
$$
f_m(k,a,b) = \sum_{p=0}^m (-1)^p c_{m-p}(a) h_p(y_{b+1}, \dots , 
y_{b+k}), \tag 15
$$
where  $h_p(z_1, \dots , z_k)$ denotes the $p^{\text{th}}$ complete 
symmetric polynomial in  $z_1, \dots , z_k$.  Equivalently, 
$$
\sum_{m} f_m(k,a,b) t^m =  \sum_{i \geq 0} c_i(a) t^i \bigg/ 
\prod_{j= b+1}^{b+k}(1+y_jt).
$$
For $k > 0$, $a \geq 0$, $b \geq 0$, set
$$
D(k,a,b) = \det\bigl(f_{a+j-i}(k,a+k-i,b) \bigr){}_{1 \leq i,j \leq k}.
\tag 16
$$
\proclaim{Proposition 4.1}  For $w = w_0^{(N)}$,  
$$
 \frak{S}_w(c,y) = \prod_{i=1}^{l-1} D(k_{i+1},@,@, n_i,@,@, n_l-
 n_{i+1}).
 $$
 \endproclaim
 \demo{Proof}
 Note that both sides of the identity to be proved are in the $R$-module 
 $M$ described in Lemma 2.1, where 
 $R = \Bbb{Z}[y_1,\dots,y_{n_l-1}]$.  By that lemma, it therefore 
 suffices to prove the formula after specializing each $c_i(j)$ to 
 $e_i(x_1,\dots,x_j)$.  The permutation $w$ is dominant, 
 so, by \cite{M, (6.14)}, 
 $\frak{S}_w(x,y) = \prod (x_p - y_q)$, the product over all $(p,q)$ that 
 appear in the diagram $D(w)$ of $w$.  That is, 
 $$
 \frak{S}_w(x,y) = \prod_{i=1}^{l-1} \prod_{p =1}^{n_i} 
 \prod_{q=n_l-n_{i+1}+1}^{n_l-n_i} (x_p - y_q).
 $$
 We want to show that the $i^{\text{th}}$ term in the product of the 
 proposition specializes to the $i^{\text{th}}$ term of this product.  
 Equivalently, we must show that $D(k,n,0)$ specializes to 
 $\prod_{p=1}^n \prod_{q=1}^k (x_p - y_q)$ for any positive $n$ and $k$.  
 But this also follows from \cite{M, (6.14)}, cf\. \cite{F2, (9.6)}. 
 \enddemo

 We next want to see what happens to the Schubert polynomials 
 $\frak{S}_w(c)$ when we carry out a specialization $c \mapsto g$.  For 
 this, it is useful to have a more graphic description of this 
 specialization.  Write $x_1, \dots, x_{n+1}$ as the vertices of the 
 Dynkin diagram for $(A_n)$, and regard $g_1[1], \dots, g_n[1]$ 
 as the edges, with $g_i[1]$ connecting $x_i$ to $x_{i+1}$.  Now regard 
 $g_i[j]$ as the path starting at vertex $x_i$, moving $j$ steps to the 
 right, and ending at vertex $x_{i+j}$; in particular, this identifies 
 $g_i[0]$ with $x_i$.  With this interpretation, $c_i(k)$ is the {\it sum 
 of all products of disjoint paths that cover exactly} $i$ {\it vertices, 
 all in} $\{x_1, \dots, x_k\}$.

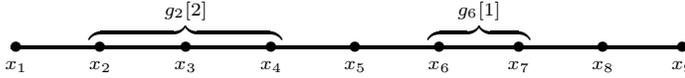
 \phantom{x}
 \vbox{$$\underset{x_1}\to\bullet\hskip .3in
 \overbrace{\underset{x_2}\to\bullet\hskip .3in
  \underset{x_3}\to\bullet\hskip .3in
   \underset{x_4}\to\bullet}^{g_2[2]} \hskip .3in
    \underset{x_5}\to\bullet\hskip .3in
    \overbrace{\underset{x_6}\to\bullet\hskip .3in
     \underset{x_7}\to\bullet}^{g_6[1]} \hskip .3in
      \underset{x_8}\to\bullet\hskip .3in
       \underset{x_9}\to\bullet\hskip .3in$$
       \vskip -.29in \moveright .57in \vbox{\hrule width 3.53in height 1pt}
       \vskip .29in}
	
	\noindent This figure illustrates, for example, that the monomial 
	$x_1g_2[2]x_5g_6[1]x_9$ appears in the expansion of $c_8(9)$.

	This definition is equivalent to the inductive definition
	$$
	c_i(k) = \sum_{j=0}^i g_{k-j}[j] \, c_{i-j-1}(k-j-1). \tag 17
	$$ 
	This inductive definition is seen to be equivalent to the matrix definition 
	given in Section 1 by expanding the determinant along the right column. 

	We consider the determinant $D(k,a,b)$ and polynomial $f_m(k,a,b)$ 
	defined before Proposition 4.1, but with each $c_i(j)$ replaced by the 
	corresponding sums of products of $g_s[t]$'s.  We need the following 
	fact.

	\proclaim{Lemma 4.2} Suppose $a$, $b$, and $k$ are given, and assume 
	that $g_i[j] = 0$ whenever $a < i+j < a+k$ and $j > 0$.  Then 
	$$
	D(k,a,b) = \det\bigl(f_{a+j-i}(k,a,b) \bigr){}_{1 \leq i,j \leq k}.
	$$ 
	\endproclaim
	\demo{Proof}
	The hypotheses imply that 
	$f_i(k,p,b) = f_i(k,p-1,b) + x_pf_{i-1}(k,p-1,b)$
	for $a < p < a+k$.  It is then a matter of elementary row reduction to 
	show that the rows of the matrix for $D(k,a,b)$ can be replaced by the 
	rows of the determinant on the right, by successively subtracting 
	linear combinations of lower rows.
	\enddemo

	We next show that the polynomials that give Giambelli formulas in 
	quantum cohomology for Schubert varieties in partial flag manifolds 
	\cite{CF2} can also be realized as specializations of the universal 
	Schubert polynomials.  With $N$ a set of $l$ positive integers as above, 
	let $F^{(N)}$ be the flag variety of flags  $V_1 \subset \dots \subset 
	V_{l-1} \subset V$, with $V = V_l$ a fixed $n_l$-dimensional vector 
	space, and $\dim(V_i) = n_i$.  For computations in quantum cohomology 
	of  $F^{(N)}$  one introduces variables $q_1, \dots , q_{l-1}$,  where 
	$q_i$ has degree $n_{i+1} - n_{i-1}$.  For a permutation $w$ in 
	$S^{(N)}$, we let $\frak{S}_w^{(N)}$ be the result of specializing in 
	$\frak{S}_w(g)$ the variables $g_i[0]$ to $x_i$, the variables 
	$$
	g_{n_{i-1}+1}[k_i+k_{i+1}-1] \quad \text {to} \quad (-1)^{k_{i+1}+1}q_i 
	\tag 18
	$$
	for $1 \leq i \leq l-1$,  and all other variables $g_i[j]$ are set equal to 
	$0$. 

	Another procedure is used in \cite{CF2}, which can be seen to amount to 
	the following.  First do the substitutions $c_i(j) \mapsto c_i(n_k)$ for 
	$j \in [n_k,n_{k+1})$, and then do the preceding substitutions from the 
	$c_i(n_k)$ to polynomials in $x$'s and $q$'s.
	 
	 \proclaim{Proposition 4.3} For $w$ in $S^{(N)}$, these two procedures 
	 give the same polynomial $\frak{S}_w^{(N)}$ in 
	 $\Bbb{Z}[x_1, \dots, x_{n_l},q_1, \dots, q_{l-1}]$.
	 \endproclaim
	 \demo{Proof}
	 When $w = w_0^{(N)}$, this follows from Proposition 4.1 and 
	 Lemma~4.2; note that the vanishing of many $g_i[j]$ guarantee that the 
	 lemma applies to each determinant 
	 $D(k_{i+1},n_i,n_l-n_{i+1})$.  The general case then follows from this 
	 case, since every other Schubert polynomial, for $w$ in $S^{(N)}$, can 
	 be obtained from the polynomial for $w_0^{(N)}$ by a sequence of 
	 operations $\partial_k^{(y)}$ acting on the $y$ variables alone, and 
	 specializing in the $c$ variables commutes with these operations.
	 \enddemo

	 There are other determinantal formulas for universal Schubert 
	 polynomials that can be deduced from corresponding formulas for 
	 classical Schubert polynomials, by the same procedure as in 
	 Proposition 4.1.   For example, this will be the case whenever all the  
	 appearances of $c_i(j)$ in entries of the matrix never have the same 
	 $j$ appearing in two different rows, or if the same $j$ never appears in 
	 two different columns.  Indeed, the determinant of such a matrix will 
	 be in the module denoted $M$ in Section 2.  In fact, as in the proof of 
	 Lemma 4.2, the same is true if the matrix can be transformed by an 
	 appropriate sequence of elementary row (resp. column) operations into 
	 a matrix of this form.

	 For example, suppose $w$ is a Grassmannian permutation with descent 
	 at $r$.  The shape of $w$ is the partition $\lambda = (w(r)-1, 
	 \dots,w(2)-2, w(1)-1)$.  Let $\mu$ = $\tilde{\lambda}$, which is the 
	 shape of $w^{-1}$.  Let $\phi$ be the flag of $w^{-1}$; this is an 
	 increasing sequence of positive integers, constructed from the Lehmer 
	 code 
	 $(c_1, \dots , c_n)$ of $w^{-1}$ by arranging the integers of the form 
	 $\max\{j \geq i \mid c_j \geq c_i\}$, taken over those $i$ with 
	 $c_i \neq 0$, in weakly increasing order (see \cite{M, p. 14}).  Let 
	 $f_m(k,a) = f_m(k,a,0)$, where $f_m(k,a,b)$ is defined in \thetag{15}.  
	 Then we have

	 \proclaim{Proposition 4.4} If $w$ is a Grassmannian permutation, then 
	 $$
	 \frak{S}_w(c,y) = \det \bigl( f_{\mu_i+j-i}(\phi_i , r+j-1) \bigr),
	 $$
	 with $\mu$ and $\phi$ as defined in the preceding paragraph. 
	 \endproclaim
	 \demo{Proof}
	 The second indices ``$r+j-1$'' occuring in the matrix guarantee that the 
	 determinant is in the $\Bbb{Z}[y]$-module $M$.  It therefore suffices to 
	 prove the corresponding formula in the classical case, i.e., to prove that 
	 $$
	 \frak{S}_w(x,y) = \det \bigl( e_{\mu_i+j-i}(\phi_i , r+j-1) \bigr),
	 $$  
	 where $e_m(u,v)$ denotes the coefficient of $t^m$ in 
	 $ \left. \prod_{i=1}^v(1+x_it) \right/ \prod_{j=1}^u(1+y_jt)$.  This 
	 formula is known (\cite{M, (6.15), (3.8)}, cf\. \cite{F2, (9.18)}, except 
	 that in these references the determinant is of the matrix 
	 $\bigl( e_{\mu_i+j-i}(\phi_i , r)\bigr)$.  But these two determinants 
	 are seen to be equal by doing elementary column operations, adding 
	 multiples of left columns to those on the right.
	 \enddemo

	 \example{Remark 4.5} When the universal Schubert polynomial 
	 $\frak{S}_w(c,y)$ is specialized to the quantum variables, one recovers 
	 the formula of Kirillov \cite{K} for Grassmannian permutations. 
	 \endexample

	 Other formulas can be proved in the same way; they give formulas for 
	 some Schubert polynomials $\frak{S}_w(c)$ as determinants of 
	 matrices with entries $c_i(j)$ where $j$ is constant in columns, but 
	 with different values in different columns.  

	 In fact, however, some experimenting indicates that there may be more 
	 determinantal formulas if one looks for such matrices with $j$ being 
	 constant in {\it rows}, with different values in different rows.  In fact, 
	 many of these polynomials have the following special form.  For 
	 arbitrary sequences $a_1, \dots, a_n$ and $b_1, \dots, b_n$ of 
	 nonnegative integers, form the $n \times n$ matrix $C(a,b)$ whose 
	 $(i,j)$ entry is $c_{a_i+j-i}(b_i)$, except that, if $a_i = 0$,  then the 
	 $(i,j)$ entry of $C(a,b)$ is the Kroneker index $\delta_{i@,@,j}$.  Note that 
	 $C(a,b)$ has the numbers  $c_{a_i}(b_i)$ appearing down the diagonal.  
	 Define $D_{a_1 \ldots  a_n}(b_1 \ldots  b_n)$ to be the determinant of 
	 this matrix $C(a,b)$.  The formulas we are looking for have the form
	 $$
	 \frak{S}_w(c) = D_{i_{\sigma(1)} \ldots i_{\sigma(n)}}(\sigma(1) 
	 \ldots \sigma(n)). \tag 19
	 $$
	 Here $w$ be a permutation in $S_{n+1}$, $I = (i_1, \dots, i_n)$ is the 
	 code $c'(w)$ defined in Section 2, and $\sigma$ is some permutation in 
	 $S_n$.  Note that the product of the diagonal terms in such a 
	 determinant is the leading term of $\frak{S}_w(c)$ found in 
	 Proposition 2.2.  Here are some examples, for permutations in $S_5$:
	 $$\alignedat 3
	 w &= 5\,1\,4\,2\,3 \qquad &I &= (1,1,2,2) \qquad &\frak{S}_w(c) &= 
	 D_{2\,2\,1\,1}(4\,3\,2\,1)\\
	 w &= 3\,5\,1\,2\,4 \qquad &I &= (0,2,2,1) \qquad &\frak{S}_w(c) &= 
	 D_{1\,0\,2\,2}(4\,1\,3\,2)\\
	 w &= 3\,2\,5\,1\,4 \qquad &I &= (1,0,3,1) \qquad &\frak{S}_w(c) &= 
	 D_{1\,0\,1\,3}(4\,2\,1\,3)
	 \endalignedat
	 $$
	 Note that the last permutation is not vexillary.  In fact, all but 8 of the 
	 120 permutations in $S_5$ have such determinantal expressions, but 
	 these do not include all of the 103 vexillary permutations. For example, 
	 $w = 1\,5\,3\,2\,4$ is vexillary, but has no such determinantal 
	 expression.  It would be interesting to characterize those permutations 
	 which have such a determinantal expression, or to give a rule for a 
	 permutation $\sigma$ (which is often not unique) for those that do.  As 
	 before, such a formula can be verified by computing its classical 
	 specialization, but the formulas seem easier to detect for the universal 
	 polynomials.

	 One can write a product of two universal Schubert polynomials as a 
	 linear combination of universal Schubert polynomials, with 
	 coefficients in the ring $\Bbb{Z}[g]$ generated by all $g_i[j]$, $j \geq 
	 0$:
	 $$
	 \frak{S}_u(g)\cdot\frak{S}_v(g) = \sum a_w(g) \frak{S}_w(g). \tag 20
	 $$
	 Even when specialized all the way to the classical case, formulas for 
	 the coefficients are only known in special cases, cf\. \cite{S}.  Many of 
	 these special cases have been extended to the quantum setting 
	 (\cite{F-G-P}, \cite{CF2}).  In these settings, the coefficients that appear 
	 with monomials in the $x$ and $q$ variables are all positive, for 
	 geometric reasons.  This is no longer the case for these universal 
	 polynomials.  Indeed, given the signs with which the $q$ variables 
	 appear in the quantum partial flag specializations, one knows that 
	 certain of these coefficients must be negative.

	 There is one special case where one does have a positive expansion, 
	 which is the case of multiplying two single terms $c_i(k)\cdot c_j(k)$.
	 Note that $c_i(k) = \frak{S}_w(c)$, where $w$ is the Grassmannian 
	 permutation in $S_{k+1}$ with descent at $k$ and $w(k+1) = k+1-i$. 

	 To state this result, let $i$, $j$, and $k$ be integers, with  $k$ positive 
	 and $0 \leq i, j \leq k$.  Let $\Cal{A}(i,j,k)$ be the set of Grassmannian 
	 permutations $w$ in $S_{k+1}$ with descent at $k-1$, such that 
	 $w(k) \leq k+1-\max(i,j)$ and $w(k)+w(k+1) = 2k+3-(i+j)$.

	 \proclaim{Proposition 4.6} For $0 \leq i, j \leq k$, 
	 $$
	 c_i(k)\cdot c_j(k) = \sum_{w \in \Cal{A}(i+1,j+1,k+1)} \frak{S}_w(c) + 
	 \sum_{w \in \Cal{A}(i,j,k)} g_k[1] \, \frak{S}_w(c) + 
	 \sum_{p=1}^{k-1} \, g_{k-p}[p+1] \, A_p,
	 $$
	 where $A_p = \sum \frak{S}_u(c)$, the sum over $u$ of the form $u = 
	 w@,@, t_{k-p,k}$,  where $w$ varies over those permutations in 
	 $\Cal{A}(i,j,k)$ for which $w(k-p) > w(k)$, and $u = w@,@, t_{k-p,k}$ is 
	 the result of interchanging the values of $k-p$ and $k$.
	 \endproclaim
	 \demo{Proof}
	 The proof is by explicit calculation, using the following facts.  First, 
	 for $w \in \Cal{A}(i,j,k)$, with $a = w(k)$ and $b = w(k+1)$,
	 $$
	 \frak{S}_w(c) = c_{k-a}(k-1)c_{k+1-b}(k) - c_{k-b}(k-1)c_{k+1-a}(k). 
	 \tag 21
	 $$
	 This is a special case of Proposition 4.4.  It then follows from (6) that 
	 for $u = w@,@, t_{k-p,k}$ as in the sum for $A_p$, 
	  $$
	  \frak{S}_u(c) = c_{k-a-p}(k-1-p)c_{k+1-b}(k) - 
	  c_{k-b-p}(k-1-p)c_{k+1-a}(k). \tag 22
	  $$
	  With these formulas, one can expand each $\frak{S}_w(c)$ and 
	  $\frak{S}_u(c)$ that occurs in the statement of the proposition.  On the 
	  other hand, from the definitions we have the formula
	  $$
	  c_s(k+1) = c_s(k) + \sum_{r=0}^{s-1} \, g_{k+1-r}[r] \, c_{s-r-1}(k-r). 
	  \tag 23
	  $$
	  One then substitutes \thetag{23} for each occurrence of any $c_s(k+1)$ 
	  in the expansion of the first sum in the formula of the proposition, and
	  verifies that one has an identity.
	  \enddemo
	   
	   \example{Remark 4.7} The first sum in the formula can be written in
	   the form
	   $$\sum_{l \geq 0}c_{i-l}(k+l)c_{j+l}(k) -
	   \sum_{l \geq 1}c_{i-l}(k)c_{j+l}(k+1).
	   $$
	   When specialized to the classical case, the other sums vanish, and the
	   proposition becomes a standard formula for multiplying elementary
	   symmetric polynomials.  \endexample
	    
	    \example{Remark 4.8} The preceding proposition can be used
	    inductively
	    -- in the explicit form given in the proof -- to write an arbitrary
	    polynomial in variables $c_i(k)$ as a linear combination of such
	    monomials so that no product $c_i(k)c_j(k)$ occurs with positive $i$ and $j$.
	    Such a polynomial can be written as a linear combination of universal 
	    Schubert
	    polynomials, using Proposition 2.2.  \endexample
	     
	     One can construct a {\it universal ring}  $R_n$, which is a natural place
	     for calculations, and which specializes to the classical and quantum
	     cohomology rings of flag and partial flag varieties.  This is an algebra
	     over the polynomial ring $\Bbb{Z}[g]$ in variables $g_i[j]$ for $j > 0$
	     and $i+j \leq n+1$.  The ring can be defined as
	     $$
	     R_n = \Bbb{Z}[g][x_1, \dots , x_{n+1}] \, / \, (c_1(n+1), \dots ,
	     c_{n+1}(n+1)), \tag 24
	     $$
	     where each $c_i(n+1)$ is the polynomial in $x$ and $g$ variables
	     defined in Section 2 or \thetag{17}.  
	     The Schubert polynomials $\frak{S}_w(g)$, as
	     $w$ varies over $S_{n+1}$, form a basis for this algebra over
	     $\Bbb{Z}[g]$.
	      
	      This ring has an obvious inner product $\langle \, , \, \rangle$, with
	      values in $\Bbb{Z}[g]$, which is obtained by multiplying and then 
	  picking off the coefficient of $x_1^nx_2^{n-1}  \dots  x_n$, or the 
	  coefficient of $\frak{S}_{w_0}(g)$.  Computations in low degrees
	  lead one to conjecture the following orthogonality:
	  $$
	  \langle \frak{S}_u(g)\, , \, \omega(\frak{S}_{v @, w_0}) \rangle =
	  \delta_{u@,v}, \tag 25
	  $$
	  where $\omega$ is the involution defined by $\omega(g_i[j]) = 
	  (-1)^{j+1}g_{n+2-i-j}[j]$.
	  Such an othogonality was difficult to prove in the quantum case 
	  (\cite{F-G-P},
	  \cite{K-M}), and promises to be even more difficult for this
	  generalization.  Similar calculations lead one to believe \footnote{A. N.
	  Kirillov reports that he can prove \thetag{25} and \thetag{26}.} that
	  $$
	  \frak{S}_w(c,c) = 0 \quad {\text { for }}\quad  w \neq 1. \tag 26
	  $$

	  In a more recent work \cite{B-F}, formulas for the general degeneracy loci
	  described in \cite{A-DF-K} and \cite{L-M} are found by different methods. 
	  When specialized to the loci considered here, they give different formulas
	  for the universal double Schubert polynomials.

	  It should be interesting to look for analogues of these polynomials for the
	  other classical groups.

\enddocument